\documentclass[12pt]{article} 
\usepackage{amssymb}
\usepackage{graphicx}
\usepackage{amsmath}
\usepackage{amsfonts}
\usepackage{hyperref}
\usepackage[utf8]{inputenc}
\DeclareUnicodeCharacter{0301}{}
\newcommand{\singlespacing}{\let\CS=\@currsize\renewcommand{\baselinestreatch}{1.0}\tiny\CS}
\newcommand{\doublespacing}{\let\CS=\@currsize\renewcommand{\baselinestreatch}{1.5}\tiny\CS}

 \newtheorem{theorem}{Theorem}[section]
 \newtheorem{lemma}{Lemma}[section]

 \DeclareMathOperator{\Lie}{\mathcal{L}}

 \numberwithin{equation}{section}
\begin{document}
\begin{center}
\Large {\textbf{On The Triviality Of $m$-Modified Conformal Vector Fields}}\\
\vskip.1in
Rahul Poddar$^1$, Ramesh Sharma$^2$\\
\vskip.1in
$^1$Harish-Chandra Research Institute, Prayagraj, 211019, UP, INDIA. E-mail:  rahul27poddar@gmail.com   \\
$^2$University of New Haven, West Haven, CT 06516, USA. E-mail: rsharma@newhaven.edu  \\

\end{center}

\begin{abstract} 
\noindent
We prove that a compact Riemannian manifold $M$ does not admit any non-trivial $m$-modified homothetic vector fields. In the corresponding case of an $m$-modified conformal vector field $V$, we establish an inequality that implies the triviality of $V$. Further, we demonstrate that an affine Killing $m$-modified conformal vector field on a non-compact Riemannian manifold $M$ must be trivial. Finally, we show that an $m$-modified gradient conformal vector field is trivial under the assumptions of polynomial volume growth and convergence to zero at infinity.
\end{abstract}
 


\section{Introduction}
In recent years, there has been an epistemic curiosity among researchers to study Ricci solitons, which correspond to self-similar solutions of the Ricci flow equation. A Riemannian manifold $(M,g)$ is said to be a Ricci soliton if it satisfies the equation 
\begin{equation}\label{1.1}
   \frac{1}{2} \Lie_V g + Ric = \lambda g,
\end{equation}
where $\Lie$ denotes the Lie-derivative operator, $V$ is a smooth vector field on $M$, $g$ is the Riemannian metric, Ric is the Ricci tensor of $g$ and $\lambda$ a real constant. If $V$ is the gradient of a smooth function $f$, i.e., $V=\nabla f$, then $(M,g,f,\lambda)$ is called a gradient Ricci soliton. As a generalization of Einstein metrics, gradient Ricci solitons, gradient Ricci almost solitons, Catino \cite{catino2012generalized} defined a generalized quasi-Einstein manifold  as an $n$-dimensional smooth Riemannian manifold $(M,g)$ satisfying the equation
\begin{equation}\label{1.2}
    Ric + Hess f - \mu  df\otimes df =\lambda g,
\end{equation}
where it was shown that a complete $(M,g)$ satisfying (\ref{1.2}) with harmonic Weyl tensor and zero radial Weyl curvature is locally a warped product with $(n-1)$-dimensional Einstein fiber. In particular, if $$\mu = \frac{1}{m}$$ for a real $m$ such that $0 < m \leq \infty$, then (\ref{1.2}) becomes   
\begin{equation}\label{1.3}
  Ric + Hess f - \frac{1}{m}  df\otimes df =\lambda g,
 \end{equation}
where the left-hand side is called the $m$-Bakry-\'Emery Ricci tensor. For $\lambda$ constant, (\ref{1.3}) defines an $m$-quasi Einstein manifold (see Case, Shu and Wei \cite{case2011rigidity}, He, Petersen and Wylie \cite{he2010classification}). 
\noindent
As a generalization of (\ref{1.3}), Limoncu \cite{limoncu} and Barros and Ribeiro, Jr. \cite{barros2012integral} considered the following equation
\begin{equation}\label{1.4}
    \frac{1}{2}\Lie_V g + Ric - \frac{1}{m}v \otimes v = \lambda g,
\end{equation}
where $v$ is the $1$-form metrically associated to the potential vector field $V$ and $\lambda$ is a constant. We note that as $m \to \infty$, (\ref{1.4}) reduces to the classical Ricci soliton equation defined by (\ref{1.1}). It has been shown by Poddar et al. \cite{poddar} that an $m$-quasi-Einstein manifold ($\lambda$ constant in (\ref{1.3})) with constant scalar curvature admitting a non-parallel closed conformal vector field is Einstein. Thus, considering the metric $g$ in (\ref{1.4}) to be both Einstein and quasi-Einstein, Zhang and Chen \cite{zhang} studied the following equation
\begin{equation}\label{1.5}
    \frac{1}{2}\Lie_V g - \frac{1}{m}v \otimes v = \sigma g,
\end{equation}
where $\sigma$ is a real constant. We will call the vector field $V$ satisfying equation (\ref{1.5}) with $0<m<\infty$, as an $m$-modified homothetic vector field. If the constant $\sigma$ is considered as a function on $M$, then the vector field $V$ is known as an $m$-modified conformal vector field, and trivial when $V=0$.  Evidently, as $m \to \infty$, (\ref{1.5}) reduces to a conformal vector field 
\begin{equation}\label{1.6}
    \Lie_V g = 2\sigma g.
\end{equation}
If $\sigma$ is constant in (\ref{1.6}), then $V$ is said to be a homothetic vector field and Killing if $\sigma = 0$. If $V=\nabla f$ in (\ref{1.5}), for some smooth function $f$ on $M$, then it is said to be an $m$-modified gradient conformal vector field and equation (\ref{1.5}) then assumes the form
\begin{equation}\label{1.7}
    Hess f - \frac{1}{m}df \otimes df = \sigma g,
\end{equation}
where $Hess$ denotes the Hessian of $f$ with respect to $g$. For more details on an $m$-modified conformal vector field, we refer the reader to the work of Zhang and Chen \cite{zhang}.\\

\noindent
The purpose of this paper is to investigate the geometry of an $m$-modified conformal vector field $V$ on a compact and non-compact smooth Riemannian manifold $M$, respectively. In the latter case, we study $V$ when it is affine Killing, followed by $V$ when it is gradient and $M$ is complete.\\

\section{$m$-Modified Conformal Vector Fields On Compact Riemannian Manifolds}
\subsection{Triviality Results}
\noindent
Zhang and Chen \cite{zhang} showed that a connected, compact Riemannian manifold of dimension at least $2$ admits no $m$-modified conformal vector field of constant length. It is also well-known that a homothetic vector field on a compact Riemannian manifold is Killing (see Yano \cite{yano1970integral}). Motivated by this result, we prove that a compact Riemannian manifold admits no non-trivial $m$-modified homothetic vector fields.
\begin{theorem}
    Any $m$-modified ($m$ finite) homothetic vector field on a compact Riemannian manifold is trivial.
\end{theorem}

\noindent
\textbf{Proof.} Tracing equation (\ref{1.5}) gives
\begin{equation}\label{2.1}
    div V = \frac{1}{m}|V|^{2} + n\sigma.
\end{equation}
Equation (\ref{1.5}) can also be exhibited as
\begin{equation}\label{2.2}
    g(\nabla_X V,Y) + g(\nabla_Y V, X) = \frac{2}{m}g(V,X)g(V,Y) + 2\sigma g(X,Y).
\end{equation}
The exterior derivative $dv$ of the $1$-form $v$ is given by
\begin{equation}\label{2.3}
    g(\nabla_X V,Y) - g(\nabla_Y V, X) = 2(dv)(X,Y).
\end{equation}
As $dv$ is skew-symmetric, we define a tensor field of type $(1,1)$ given by
\begin{equation}\label{2.4}
    dv(X,Y) = g(FX,Y).
\end{equation}
Evidently, $F$ is skew self-adjoint, i.e., $g(FX,Y)=-g(X,FY)$. Thus, (\ref{2.4}) assumes the form
\begin{equation}\label{2.5}
    g(\nabla_X V,Y) - g(\nabla_Y V, X) = 2g(FX,Y).
\end{equation}
Adding equations (\ref{2.2}) and (\ref{2.5}) side by side, and factoring out $Y$ gives
\begin{equation}\label{2.6}
    \nabla_X V = \frac{1}{m}g(V,X)V + \sigma X + FX.
\end{equation}
Now, we compute the covariant derivative of the squared $g$-norm of $V$ using (\ref{2.6}) as follows
\begin{equation}\label{2.7}
    \nabla_X |V|^{2} = 2g(\nabla_X V, V) = \frac{2}{m}g(V,X)|V|^{2} + 2\sigma g(X,V) + 2g(FX,V).
\end{equation}
Factoring out $X$ from (\ref{2.7}) gives
\begin{equation*}
    \nabla|V|^{2} = \frac{2}{m}|V|^{2}V + 2\sigma V - 2FV.
\end{equation*}
Taking the divergence on both sides of the preceding equation entails
\begin{equation*}
    \Delta |V|^{2} = \frac{2}{m}div(|V|^{2}V) + 2\sigma div V - 2div FV.
\end{equation*}
i.e.,
\begin{equation*}
    \frac{1}{2} \Delta |V|^{2} = \frac{1}{m}(V|V|^{2} + |V|^{2}div V) + \sigma div V - div FV.
\end{equation*}
Substituting $X$ with $V$ in (\ref{2.7}), using it in the above equation, and further simplification leads us to
\begin{equation}\label{2.8}
\frac{1}{2} \Delta |V|^{2} = \bigg(\frac{2}{m}|V|^{2} + div V \bigg)\bigg(\frac{1}{m}|V|^{2} + \sigma \bigg) - 2div FV.
\end{equation}
Eliminating $|V|^{2}$ between (\ref{2.1}) and (\ref{2.8}) shows that
\begin{equation}\label{2.9}
    \frac{1}{2} \Delta |V|^{2} = 3(div V)^{2} + (3-5n)\sigma div V + 2n(n-1)\sigma^{2} - 2div FV.
\end{equation}
Integrating equation (\ref{2.9}) over compact $M$ and using divergence theorem provides, $div V = 0$ and $\sigma =0$. The use of the foregoing equations in (\ref{2.1}) at once gives $V=0$, i.e., $V$ is trivial. This completes the proof.\\

\noindent
In order to set the context for our next result, we recall the following result \cite{yano1952harmonic} of Yano ``If a compact Riemannian manifold $(M,g)$ admits a conformal vector field $V$ satisfying $Ric(V,V) \leq 0$, then $V$ is parallel, and if the Ricci curvature is negative definite, then the conformal vector field $V$ is zero". Motivated by this and Theorem 2.1, we examine the case when $V$ is an $m$-modified conformal vector field, and establish the following triviality result.

\begin{theorem}
    Any $m$-modified ($m$ finite) conformal vector field $V$ on a compact Riemannian manifold satisfying the inequality
    \begin{equation*}
        Ric(V,V) \leq \frac{2(1-n)}{nm^2}(max |V|)^{4}
    \end{equation*}
is zero, and hence trivial.
\end{theorem}

\noindent
\textbf{Proof.} Let us recall the following integral formula (formula 1.11, Yano \cite{yano1952harmonic}) for any smooth vector field $V$ on a compact Riemannian manifold $(M,g)$:
\begin{equation}\label{2.10}
    \int_M [Ric(V,V) + \frac{1}{2}|\Lie_V g|^{2} - |\nabla V|^{2} - (div V)^{2}]dv_g = 0,
\end{equation}
where $dv_g$ is the volume element of $(M,g)$. Tracing the equation
\begin{equation}\label{2.11}
     \Lie_V g  = 2\sigma g + \frac{2}{m}v \otimes v,
\end{equation}
of an $m$-modified conformal vector field $V$, we get
\begin{equation}\label{2.12}
    div V = \frac{1}{m}|V|^{2} + n\sigma.
\end{equation}
Next, using (\ref{2.11}) we compute
\begin{equation}\label{2.13}
    |\Lie_V g|^{2} = 4n\sigma^{2} + \frac{4}{m^2}|V|^{4} + \frac{8\sigma}{m}|V|^{2}.
\end{equation}
Squaring (\ref{2.12}) provides
\begin{equation*}
    \frac{2}{n}(div V)^{2} = 2n\sigma^{2} + \frac{2}{nm^2}|V|^{4} + \frac{4\sigma}{m}|V|^{2}.
\end{equation*}
The use of this equation in (\ref{2.13}) yields 
\begin{equation*}
    \frac{1}{2}|\Lie_V g|^{2} = \frac{2}{n}(div V)^{2} + \frac{2}{m^2}\bigg(1-\frac{1}{n}\bigg)|V|^{4}.
\end{equation*}
Substituting this value of $\frac{1}{2}|\Lie_V g|^{2}$ in (\ref{2.10}), we get
\begin{equation*}
    \int_M [Ric(V,V) + \frac{2}{m^2}\bigg(1-\frac{1}{n}\bigg)|V|^{4}]dv_g = \int_M [|\nabla V|^{2} + \bigg(1 - \frac{2}{n}\bigg)(div V)^{2}]dv_g.
\end{equation*}
This implies that
\begin{equation}\label{2.14}
    \int_M Ric(V,V) dv_g \geq \int_M \frac{2}{m^2}\bigg(\frac{1}{n} - 1 \bigg)|V|^{4}dv_g.
\end{equation}
Noting that $|V|^{4} \leq (max |V|)^{4}$, the above inequality implies that
\begin{equation*}
    \int_M Ric(V,V) dv_g \geq \int_M \frac{2}{m^2}\bigg(\frac{1}{n} - 1 \bigg)(max |V|)^{4}dv_g.
\end{equation*}
At this point, using our hypothesis, we conclude that $max |V| = 0$ and so $V=0$, completing the proof.\\

\section{$m$-Modified Conformal Vector fields On Non-compact Riemannian Manifolds}

\subsection{Triviality Results}
\subsubsection{Triviality Of An Affine Killing $m$-modified Conformal Vector Field}
\noindent
 We know that a conformal vector field $V$ which is also affine Killing (i.e., satisfies $\Lie_V \nabla = 0$, and preserves the geodesics along with their affine parameters), is homothetic. This can be easily verified by the integrability condition (Yano \cite{yano1970integral}):
\begin{equation}\label{3.1}
    (\Lie_V \nabla)(X,Y) = (X\sigma)Y + (Y\sigma)X - g(X,Y)\nabla \sigma,
\end{equation}
for a conformal vector field $V$: $\Lie_V g = 2\sigma g$. Using $\Lie_V \nabla = 0$ in (\ref{3.1}) and contracting with respect to $X$ and $Y$ shows that $\sigma$ is constant, i.e., $V$ is homothetic. This suggests the following question: ``What can be said about an $m$-modified ($m$ finite) conformal vector field $V$ which is also affine Killing?". It turns out that this condition places a severe restriction on $V$ so much so that $V$ becomes zero, i.e., trivial. More precisely, we prove
\begin{theorem}
    An $m$-modified ($m$ finite) conformal vector field which is also affine Killing, is trivial.
\end{theorem}

\noindent
\textbf{Proof.} First, tracing the $m$-modified conformal vector field equation
\begin{equation}\label{3.2}
     \Lie_V g  = 2\sigma g + \frac{2}{m}v \otimes v,
\end{equation}
we get
\begin{equation}\label{3.3}
    div V = n\sigma + \frac{1}{m}|V|^{2}.
\end{equation}
Next, we recall the following commutation formula (Yano \cite{yano1970integral}):
\begin{eqnarray*}
&&(\Lie_V\nabla_X g - \nabla_X \Lie_V g - \nabla_{[V,X]}g)(Y,Z) \nonumber\\ &=& -g((\Lie_V\nabla)(X,Y), Z) - g((\Lie_V\nabla)(X,Z), Y)\nonumber
\end{eqnarray*}
$\forall \: V,X,Y,Z \in\mathfrak{X}(M)$. Noting that $\nabla g = 0$, and $\Lie_V \nabla = 0$ by hypothesis, we see that the commutation formula implies that $\nabla_X \Lie_V g = 0$, i.e., $\Lie_V g$ is parallel. Using (\ref{3.2}) in it, we get
\begin{equation}\label{3.4}
    (X\sigma)g(Y,Z) + \frac{1}{m}[(\nabla_X v)(Y)v(Z) + v(Y)(\nabla_X v)Z] = 0.
\end{equation}
Contracting it with respect to $Y$ and $Z$ gives 
\begin{equation*}
    n X\sigma + \frac{1}{m}\nabla_X|V|^{2} = 0
\end{equation*}
and hence
\begin{equation}\label{3.5}
    n\sigma + \frac{1}{m}|V|^{2} = c \:(a\: constant)
\end{equation}
Contracting (\ref{3.4}) at $X$ and $Z$ gives
\begin{equation}\label{3.6}
    Y\sigma + \frac{1}{m}[(\nabla_V v)Y + (div V)v(Y)] = 0.
\end{equation}
Equations (\ref{3.3}) and (\ref{3.5}) imply that $div V = c$. Using this in (\ref{3.6}) gives
\begin{equation}\label{3.7}
    \nabla \sigma + \frac{1}{m}[\nabla_V V + cV] = 0.
\end{equation}
Now equation (\ref{3.2}) operated on $X$, $V$ ($X$ arbitrary) provides
\begin{equation*}
    \frac{1}{2}\nabla|V|^{2} + \nabla_V V = 2\bigg(\sigma + \frac{1}{m}|V|^{2}\bigg)V.
\end{equation*}
Eliminating $\nabla_V V$ between this equation and (\ref{3.7}), we obtain
\begin{equation*}
    m(n+2)\nabla \sigma = 2[2(n-1)\sigma - 3c]V.
\end{equation*}
Its inner product with $V$ gives
\begin{equation}\label{3.8}
    m(n+2)V\sigma = 2[2(n-1)\sigma - 3c]|V|^{2}.
\end{equation}
Now, from equation (\ref{3.2}) we get
\begin{equation*}
    2g(\nabla_V V) = 2\sigma |V|^{2} + \frac{2}{m}|V|^{4}
\end{equation*}
i.e.,
\begin{equation*}
    V(|V|^{2}) = 2|V|^{2}\bigg(\sigma + \frac{1}{m}|V|^{2}\bigg)
\end{equation*}
Using (\ref{3.5}) in the above, we find that
\begin{equation*}
    n V\sigma = 2(n\sigma - c)[(1-n)\sigma + c] 
\end{equation*}
Eliminating $V\sigma$ between this and (\ref{3.8}) we obtain
\begin{equation*}
    (n\sigma - c)\bigg[\frac{n+2}{n}((1-n)\sigma + c) - ((1-n)\sigma + 3c)\bigg] = 0.
\end{equation*}
So, either $\sigma = \frac{c}{n}$, or $- \frac{c}{n-2}$. That is $\sigma$ is constant, and hence from (\ref{3.5}), $|V|$ is constant. Turning our attention back to equation (\ref{3.4}), we see that it reduces to 
\begin{equation*}
    (\nabla_X v)(Y)V + v(Y)\nabla_X V = 0.
\end{equation*}
Its inner product with $V$ yields $(\nabla_X v)Y = 0$ or $V=0$. In the first case, $\nabla_X V = 0$ which implies $\Lie_V g = 0$ and hence $div V = 0$ and therefore $c=0$. This, combined with (\ref{3.5}) shows
\begin{equation}\label{3.9}
    n\sigma  + \frac{1}{m}|V|^{2} = 0.
\end{equation}
Also, equation (\ref{3.2}) gives
\begin{equation*}
    \sigma |V|^{2} + \frac{1}{m}|V|^{4} = 0
\end{equation*}
i.e.,
\begin{equation*}
    |V|^{2}\bigg[\sigma + \frac{1}{m}|V|^{2}\bigg] = 0.
\end{equation*}
As $V$ is parallel, $|V|$ is constant. Hence, from the above equation we conclude that, either $V=0$, or $\sigma  + \frac{1}{m}|V|^{2} = 0$. The last case, in conjunction with (\ref{3.9}) shows that $\sigma = 0$ and $V = 0$. Thus, $V$ is trivial, completing the proof.\\

\subsubsection{A Structure Equation For An $m$-Modified Gradient Conformal Vector Field}
\noindent
For $0<m<\infty$, consider $u = e^{-\frac{f}{m}}$. Then we have,
\begin{equation*}
    \nabla u = -\frac{1}{m}e^{-\frac{f}{m}}\nabla f,
\end{equation*}
i.e.,
\begin{equation*}
    \frac{m}{u} Hess u = \frac{1}{m}df \otimes df - Hess f.
\end{equation*}
Therefore the fundamental equation of an $m$-modified gradient conformal vector field, given by (\ref{1.7}), can be written as 
\begin{equation}\label{3.10}
    Hess u = \lambda g
\end{equation}
where $\lambda = -\frac{1}{m}e^{-\frac{f}{m}}\sigma$ is a smooth function on $M$. Thus, in the gradient case, when $m$ is finite, we will use equation (\ref{3.10}) to study (\ref{1.7}).\\

\noindent
\textbf{Remark.} A compact Riemannian manifold with constant scalar curvature admitting a non-Killing gradient conformal vector field is isometric to a sphere (see Yano \cite{yano1970integral}). In view of equation (\ref{3.10}), this result of Yano is also true for any compact Riemannian manifold with constant scalar curvature admitting an $m$-modified gradient conformal vector field with $\sigma \neq 0$.\\

\noindent
We now derive a structure equation for an $m$-modified gradient conformal vector field with associated function $u = e^{-\frac{f}{m}}$ on a Riemannian manifold $M$. 
\begin{lemma}
    Let $(M,g)$ be a Riemannian manifold admitting an $m$-modified gradient conformal vector field. Then
    \begin{equation}\label{3.11}
         \frac{1}{2}\Delta |\nabla u|^{2} = |Hess u|^{2} - \frac{1}{n-1}Ric(\nabla u, \nabla u)
    \end{equation}
    where $u = e^{-\frac{f}{m}}$.
\end{lemma}

\noindent
\textbf{Proof.} Let us express equation (\ref{3.10}) as 
\begin{equation*}
    \nabla_X \nabla u = \lambda X,
\end{equation*}
for an arbitrary smooth vector field $X$.
Using the above equation and the formula $R(Y,X)Z = \nabla_Y\nabla_X Z - \nabla_X\nabla_Y Z - \nabla_{[Y,X]}Z$, we compute $R(Y,X)\nabla u$, to get
\begin{equation}\label{3.12}
    R(Y,X)\nabla u = (Y\lambda)X - (X\lambda)Y.
\end{equation}
Contracting (\ref{3.12}) with respect to $Y$ gives 
\begin{equation*}
    Ric(X,\nabla u) = (1-n)g(\nabla \lambda, X).
\end{equation*}
Substituting $\nabla u$ for $X$, in the above equation gives
\begin{equation}\label{3.13}
    Ric(\nabla u,\nabla u) = (1-n)g(\nabla \lambda, \nabla u).
\end{equation}
Now, we take the $g$-trace of (\ref{3.10}) to get $\Delta u = n\lambda$, and use it in the well-known Bochner formula:
\begin{equation*}
     \frac{1}{2}\Delta |\nabla u|^{2} = |Hess u|^{2} + Ric(\nabla u, \nabla u) + g(\nabla \Delta u, \nabla u)
\end{equation*}
to get,
\begin{equation}\label{3.14}
     \frac{1}{2}\Delta |\nabla u|^{2} = |Hess u|^{2} + Ric(\nabla u, \nabla u) + ng(\nabla \lambda, \nabla u)
\end{equation}
Comparing equations (\ref{3.13}) and (\ref{3.14}), we obtain
\begin{equation*}
         \frac{1}{2}\Delta |\nabla u|^{2} = |Hess u|^{2} - \frac{1}{n-1}Ric(\nabla u, \nabla u).
\end{equation*}
This completes the proof.\\

\subsubsection{Triviality Of $m$-Modified Gradient Conformal Vector Field Via Polynomial Volume Growth}

\noindent
Given a polynomial function $\sigma: (0, \infty) \to (0, \infty)$, we say that $M$ has polynomial volume growth like $\sigma(r)$ if there exists $p \in M$ such that 
\begin{equation*}
    Vol(B(p,r)) = \mathcal{O}(\sigma(r)),
\end{equation*}
as $r \to \infty$, where $Vol$ denotes the volume related to the metric $g$ and $B(p,r)$ is the geodesic ball centered at $p$ with radius $r$. Our next result pertains to the triviality of an $m$-modified gradient conformal vector field whenever $M$ has polynomial volume growth under certain Ricci curvature condition. More precisely, we establish the following

\begin{theorem}
    Let $(M, g)$ be a connected, oriented, complete, non-compact Riemannian manifold admitting an $m$-modified gradient conformal vector field and whose Ricci curvature satisfies $Ric \leq -\alpha g$, for some positive constant $\alpha$. If $M$ has polynomial volume growth and $|\nabla u|$, $|\nabla \nabla u| \in L^{\infty}(M)$, where $u = e^{-\frac{f}{m}}$, then $f$ is constant (trivial).
\end{theorem}

\noindent
\textbf{Proof.} We prove it by contradiction. Suppose $f$ is non-constant on $M$. Thus, $u$ is non-constant on $M$. Let us consider the function $w=|\nabla u|^{2}$ and a smooth vector field $X = \nabla |\nabla u|^{2}$. In view of (\ref{3.11}) and our hypothesis, we have 
\begin{equation*}
    div X = \Delta |\nabla u|^{2} = 2\bigg\{|Hess u|^{2} - \frac{1}{n-1}Ric(\nabla u, \nabla u)\bigg \} \geq \frac{2\alpha}{n-1}w.
\end{equation*}
Also,
\begin{equation*}
    g(\nabla w, X) = |\nabla|\nabla u|^{2}|^{2} \geq 0.
\end{equation*}
Now, from Kato's inequality, we get
\begin{equation*}
    |X|= 2|\nabla u||\nabla|\nabla u|| \leq 2|\nabla u||\nabla \nabla u| < \infty,
\end{equation*}
since $|\nabla u|$, $|\nabla \nabla u| \in L^{\infty}(M)$. As $M$ has polynomial volume growth, we invoke the following result of  Al\'{i}as et al. \cite{alias2021maximum}, ``Let $M$ be a connected, oriented, complete, non-compact Riemannian manifold and $X$ be a bounded smooth vector field on $M$. Let $w$ be a non-negative smooth function on $M$ such that $g(\nabla w, X)\geq 0$ and $div X \geq \alpha' w$ on $M$, for some positive constant $\alpha'$. If $M$ has polynomial volume growth, then $w$ vanishes identically on $M$", to conclude that $|\nabla u|^{2}$ vanishes identically on $M$, thereby arriving at a contradiction. Therefore, $u$ is constant and consequently $f$ is constant. This completes the proof.\\

\subsubsection{Triviality Of $m$-Modified Gradient Conformal Vector Field Via Convergence to Zero At Infinity}

\noindent
 A continuous function $u \in C^{0} (M)$ is said to converge to zero at infinity if it satisfies the condition
\begin{equation*}
    \lim_{d(x, x_0) \to \infty} u(x) = 0,
\end{equation*}
where $d(., x_0):M \to [0, \infty)$ denotes the Riemannian distance of a complete, non-compact Riemannian manifold $M$ measured from a fixed point $x_0 \in M$. In the following theorem, we prove the triviality of an $m$-modified gradient conformal vector field by way of convergence to zero at infinity of the function $|\nabla u|$, for $u=e^{-\frac{f}{m}}$, under certain Ricci curvature condition.\\

\noindent
\begin{theorem}
Let $(M, g)$ be a complete, non-compact, Riemannian manifold admitting an $m$-modified gradient conformal vector field and $Ric(\nabla u, \nabla u) \leq 0$ and $|\nabla u|$ converges to zero at infinity, where $u = e^{-\frac{f}{m}}$. Then, $f$ is constant (trivial).\\
\end{theorem}

\noindent
\textbf{Proof.} We prove this by contradiction. Suppose $f$ is not a constant on $M$. Then, $u$ is a non-constant function on $M$. Consider the function $w = |\nabla u|^{2}$, which is non-negative, non-identically vanishing, and converges to zero at infinity. Now, consider the smooth vector field $X= \nabla |\nabla u|^{2}$ on $M$. Thus, we have
\begin{equation*}
    g(\nabla w, X) = |\nabla|\nabla u|^{2}|^{2} \geq 0.
\end{equation*}
From equation (\ref{3.11}) and our hypothesis that $Ric(\nabla u, \nabla u) \leq 0$, we observe that the right-hand side of (\ref{3.11}) is non-negative.
Thus,
\begin{equation*}
    div X = \Delta |\nabla u|^{2} = 2\bigg\{|Hess u|^{2} - \frac{1}{n-1}Ric(\nabla u, \nabla u)\bigg \} \geq 0.
\end{equation*}
Hence, invoking the result of Al\'{i}as et al. \cite{alias2019maximum}, ``Let $(M,g)$ be a complete, non-compact Riemannian manifold and $X$ be an arbitrary smooth vector field on $M$. Assume that there exists a non-negative, non-identically vanishing function $w \in C^{\infty} (M)$ which converges to zero at infinity and $g(\nabla w, X) \geq 0$. If $div X \geq 0$ on $M$, then $g(\nabla w, X) \equiv 0$ on $M$", we get, $w=|\nabla u|^{2}$ is a constant. Since it converges to zero at infinity, $w\equiv 0$ on $M$, and hence $f$ is constant on $M$, a contradiction to our initial assumption regarding $f$. Therefore, $f$ is constant. This completes the proof.\\

\section{Acknowledgements} 
The first author was funded by the Department of Atomic Energy (DAE), Government of India, in the form of a Postdoctoral Research Fellowship.

\section{Statements and Declarations}

\textbf{Competing Interests}: On behalf of all the authors, the corresponding author states that there is no conflict of interest. \\

\noindent
\textbf{Data Availability}: Data sharing is not applicable to this article as no datasets were generated or analysed during the current study.\\

\end{document}